\documentclass[12pt]{article}
\usepackage{amssymb}
\usepackage{latexsym}
\def\part#1{\frac{\partial\phantom{q}}{\partial#1}}

\newenvironment{rmk}{\begin{trivlist}\item[]{\bf Remark:} }
{\end{trivlist}}

\newenvironment{prf}{\begin{trivlist}\item[]{\bf Proof:} }
{\hfill $\Box$ \end{trivlist}}

\newtheorem{thm}{Theorem}

\newtheorem{prp}[thm]{Proposition}

\newcommand{\CP}{{\mathbf C}{\rm P}}

\begin{document}
\title{Bihermitian metrics on Del Pezzo surfaces}
 \author{Nigel Hitchin\\[5pt]}
\maketitle
\section{Introduction}
Algebraic geometry and symplectic geometry are usually linked through the study of K\"ahler metrics. In this paper we shall describe an example of another connection, which is a single subject with two realizations -- bihermitian geometry \cite{Ap} and generalized K\"ahler geometry \cite{Gu1}. We shall follow here the first approach.

A bihermitian structure on a manifold $M$ for us consists of a pair of integrable complex structures $I^+$ and $I^-$, a Riemannian metric $g$ which is hermitian with respect to both, and a closed 3-form $H$ such that $\nabla^+I^+=0=\nabla^-I^-$ where $\nabla^{\pm}$ is the Riemannian connection with skew torsion $\pm H$. This geometry first appeared in the physics paper \cite{Ro} in the context of the $(2,2)$ supersymmetric sigma model.

For a K\"ahler manifold with complex structure $I$, taking $I^+=I, I^-=\pm I$ and $H=0$ solves these equations, but for some time other compact examples were scarce. This has since changed  somewhat thanks to \cite{H1}, \cite{GC} and \cite{T}. Here we shall give a concrete construction which differs from these examples and holds for any Del Pezzo surface.

 Recall that a Del Pezzo surface is defined as an algebraic surface with ample anti-canonical bundle $K^*$. Ampleness means that the line bundle $K^*$ has a hermitian metric whose curvature form $F$ is positive, defining a K\"ahler metric $g_0$. This is the data for our construction -- we take a holomorphic section $\sigma$ of $K^*$, and the function $f=\log\Vert \sigma \Vert^2$. This function is singular on the zero set  of $\sigma$ (which is an elliptic curve). The section $\sigma$ can also be seen as a holomorphic Poisson structure on $M$ and its real part is a real Poisson structure. We use this Poisson structure to define a Hamiltonian vector field for $f$, which  turns out to be well-defined and smooth on the whole of $M$. Integrating it for a time $t$ we get a Poisson diffeomorphism $\varphi_t$ and   take $I^+=I$, the original complex structure on the Del Pezzo, and $I^-=\varphi_t^*I$. For small enough $t$ we show how to define  a bihermitian metric $g$ canonically from the {\it imaginary} part of $\sigma$.

The deformation parameter $t$ has an intriguing role: clearly as $t\rightarrow 0,$ $I^-(t)\rightarrow I$, but also the metric $g(t)/t$ tends to the K\"ahler metric $g_0$. More interestingly, $\varphi_t$ preserves the elliptic curve and acts as a translation proportional to $t$. This data coincides with the algebraic version of noncommutative geometry espoused by Artin, van den Bergh and Stafford (see \cite{Staff}), so that our deformation of a K\"ahler metric to a bihermitian one tracks a deformation from a commutative Del Pezzo surface to a noncommutative one. The possibility of a link between these two areas was first conjectured by Gualtieri, and our construction provides an example of this phenomenon. Further work is clearly needed to understand the relation between the differential and algebraic geometry.
\vskip .5cm
The author wishes to thank V. Apostolov, G. Cavalcanti and M. Gualtieri for useful conversations and EPSRC for support.

\section{Bihermitian 4-manifolds}

We assemble here briefly the essential facts about bihermitian geometry. Since we are working on a 4-dimensional manifold $M$, the paper \cite{Ap} by Apostolov et al. is the best reference; \cite{H1} includes similar material for arbitrary dimensions, and is more motivated by the generalized geometry approach. The definition of a bihermitian structure in \cite{Ap} is 
 a pair of integrable complex structures $I^+$ and $I^-$ on a 4-manifold defining the same orientation and a Riemannian metric $g$ which is hermitian with respect to both. When $M$ is compact with even first Betti number, and $\omega^+,\omega^-$ are the two hermitian forms, the authors show that $d^c_-\omega^-= H =-d^c_+\omega^+$ for a closed 3-form $H$ which is equivalent (see \cite{Gu1}) to the definition in the introduction.

Given a bihermitian metric, the fundamental object to study is the 2-form
\begin{equation}
\phi(X,Y)=g([I^+,I^-]X,Y)
\label{f}
\end{equation}
A simple calculation shows that $\phi(I^+X,I^+Y)=-\phi(X,Y)=\phi(I^-X,I^-Y)$. This means that, with respect to either complex structure, $\phi$ is of type $(2,0)+(0,2)$. Thus if we define $\phi'(X,Y)=\phi(I^+X,Y)$ then $\phi+i\phi'$ is a $(0,2)$-form with respect to $I^+$. Using the hermitian metric, this section of $\Lambda^2 \bar T^*$ can be identified with a section $\sigma$ of $\Lambda^2T = K^*$ and it is shown in \cite{Ap} that this is {\it holomorphic}. A similar situation holds of course for the complex structure $I^-$. 

It follows that $\phi$ vanishes on an anticanonical divisor, and this leads as in \cite{Ap} to an appeal to the classification of surfaces for candidates to admit bihermitian metrics.

 A Del Pezzo surface has a concrete description as an algebraic surface; it is either $\CP^2$ or $\CP^1\times \CP^1$ or is obtained by blowing up $k\le 8$ points in $\CP^2$ such that no three are collinear and no six lie on a conic. By Riemann-Roch we have for the latter $\dim H^0(M,K^*)=10-k$ and $\dim H^0(\CP^1\times \CP^1,K^*)=9$, so we always have  anticanonical sections $\sigma$. In fact they are very concrete -- in the first case the divisor is the proper transform of a plane cubic curve that passes through the $k$ points.

It is convenient to work with the meromorphic section $\sigma^{-1}$ of $K$. This is a meromorphic 2-form $\omega+i\omega'$ with a pole on the anticanonical divisor. It is related to $\phi$ by:

\begin{prp} \label{P1} $$\omega+i\omega'=\frac{1}{2\Vert \phi \Vert^2} (\phi-i\phi').$$
\end{prp}
\begin{prf} In local holomorphic coordinates $\sigma =s \partial/\partial z_1\wedge \partial/\partial z_2$ and then
$$\phi+i\phi'=s (\det h_{i\bar j}) d\bar z_1\wedge d\bar z_2.$$
But
$\Vert \sigma \Vert^2= s\bar s \det h_{i\bar j}$ so
$$\phi-i\phi'= s^{-1} \Vert \sigma \Vert ^2 dz_1\wedge dz_2=2\Vert \phi \Vert^2 (\omega+i\omega').$$
\end{prf}

Where $\sigma$ is non-vanishing we therefore have a closed form $\omega+i\omega'$ from the complex structure $I^+$ and  $\omega+i\omega''$ from the complex structure $I^-$. Since each is of type $(2,0)$ with respect to its own complex structure, we have $(\omega+i\omega')^2=0=(\omega+i\omega'')^2$ or 
$$\omega^2=\omega'^2=\omega''^2,\qquad \omega\omega'=0=\omega\omega''.$$ Conversely, closed forms satisfying these constraints define a pair of integrable complex structures (we simply define the $(1,0)$ forms for $I^+$ to be those annihilated by $\omega+i\omega'$). Less obviously, the hermitian metric is defined by this data. This is the content of Theorem 2 in \cite{Ap} but more concretely we have:

\begin{prp} \label{P2} For a  bihermitian metric, the hermitian form $\omega^+(X,Y)=g(I^+X,Y)$ is the $(1,1)$ component, with respect to $I^+$,  of $4\omega''$.
\end{prp}
\begin{prf} From (\ref{P1}) $\omega''= -\phi''/2\Vert \phi\Vert^2$, and 
$$\phi''(X,Y)=\phi(I^-X,Y)=g([I^+,I^-]I^-X,Y).$$
Now $I^+$ and $I^-$ define, using the metric,  self-dual 2-forms on $M$ and so their action on the tangent space is through the Lie algebra of $SU(2)$, or the imaginary quaternions. For two imaginary quaternions $u,v$ we have
$$uv+vu=-2(u,v)1,\qquad \vert uv-vu\vert^2=4(\vert u\vert^2\vert v\vert^2-(u,v)^2).$$ 
Thus 
$$\vert\phi \vert^2=\vert[I^+,I^-]\vert^2=4(1-(I^+,I^-)^2)=4(1-p^2)$$
and 
$$I^+I^-+I^-I^+=-2p1.$$
Using this latter relation, 
$$g([I^+,I^-]I^-X,Y)=-g(I^+X,Y)-g(I^-I^+I^-X,Y)=-2g(I^+X,Y)+2pg(I^-X,Y)$$
and so
\begin{equation}
\phi''=-2\omega^++2p\omega^-
\label{fi1}
\end{equation}
Similarly
\begin{equation}
\phi'=2\omega^--2p\omega^+
\label{fi2}
\end{equation}
Eliminating $\omega^-$ gives 
$$\phi''=2(p^2-1)\omega^++p\phi'$$
and since $\phi'$ is of type $(2,0)+(0,2)$ with respect to $I^+$, the $(1,1)$ part of $\phi''$ is $2(p^2-1)\omega^+$. 

From (\ref{P1}) 
$$\omega'' = -\frac{\phi''}{2\Vert \phi\Vert^2} = -\frac{\phi''}{8(1-p^2)}$$
so that
$$(\omega'')^{1,1}=\omega^+/4.$$
\end{prf}

\section{Potentials}
In the last section we have seen that a bihermitian metric is determined in a concrete fashion from the two closed forms $\omega+i\omega'$ and $\omega+i\omega''$, holomorphic symplectic with respect to $I^+$ and $I^-$ respectively. Both real and imaginary parts are real symplectic forms. 

Locally, of course, any two complex structures are equivalent under a diffeomorphism. Moreover, the holomorphic version of the Darboux theorem says that  any two holomorphic symplectic forms are equivalent. it follows that there is a local diffeomorphism $\varphi$ such that 
$\varphi^*(\omega+i\omega')=\omega+i\omega''$, or
$$\varphi^*\omega=\omega,\qquad \varphi^*\omega'=\omega''.$$
Conversely, given a holomorphic symplectic form $\omega + i\omega'$ and a diffeomorphism $\varphi$ which is symplectic with respect to $\omega$, we can  define $\omega''=\varphi^*\omega'$ and get a bihermitian metric, so long as $(\omega'')^{1,1}$ is a {\it positive} hermitian form. In particular, given a smooth function $f$, we can form its Hamiltonian vector field $X$ and generate a local one-parameter group of symplectic diffeomorphisms $\varphi_t$.

\begin{prp} If $f$ is a K\"ahler potential on an open set $U$ then, on a compact subset,  the hermitian form  $(\varphi_t^*\omega')^{1,1}$ is positive for sufficiently small $t$.
\end{prp}
\begin{prf} Differentiating with respect to $t$ at $t=0$ gives 
$$\frac{\partial }{\partial t}\varphi_t^*\omega'\vert_{t=0}={\mathcal L}_X\omega'=d(i_{X}\omega').$$

But $(\omega+i\omega')(X,Y)=i(\omega+i\omega')(X,I^+Y)$ and $i_{X}\omega=df$ 
 so that 
$i_{X}\omega'=I^+df$. Hence 
$$\frac{\partial }{\partial t}\varphi_t^*\omega'\vert_{t=0}=dI^+df=dd^cf.$$
The function $f$ is a K\"ahler potential if and only if $dd^cf$ is positive, which means that 
$$\frac{\partial }{\partial t}(\varphi_t^*\omega')^{1,1}\vert_{t=0}$$
is positive. At $t=0$, $\varphi_t$ is the identity and since $\omega'$ is of type $(2,0)+(0,2)$, we have $0=(\omega')^{1,1}=(\varphi_0^*\omega')^{1,1}$. It follows that  $(\varphi_t^*\omega')^{1,1}$ is positive for small $t$.
\end{prf}
 We see here that, if we start with a choice of holomorphic symplectic form defining the holomorphic Poisson structure, the data for a local bihermitian structure appears to be  the same as for a K\"ahler structure. There is a difference however: whereas two K\"ahler potentials $f,f'$ give the same metric if they {\it differ} by the real part of a holomorphic function, the bihermitian structures are the same if the {\it composition} $\varphi_t\circ \varphi'_{-t}$ is holomorphic. Bihermitian geometry is a  nonlinear form  of K\"ahler geometry.
\begin{rmk} The use of Hamiltonian functions to generate  bihermitian metrics in this fashion began with an observation of D. Joyce in \cite{Ap}. A more extensive version in arbitrary dimensions is given in \cite{Ro2}.
\end{rmk}

\section{Del Pezzo surfaces}

We shall now use the local construction in the previous section on a  large open set in a Del Pezzo surface. We choose a hermitian metric on the holomorphic anticanonical line bundle $K^*$ and a holomorphic section $\sigma$. By the adjunction formula $\sigma$ vanishes on an elliptic curve $C$ (which may be degenerate). We take the open set $U=M\setminus C$ and the function $f$ on $U$ given by 
$$f=\log \Vert \sigma \Vert^2.$$
We have a holomorphic symplectic form $\sigma^{-1}=\omega+i\omega'$ on $U$ and we consider the Hamiltonian vector field $X$ of $f$ with respect to $\omega$. We shall show that $X$ is well-defined on the whole of $M$.

Write $X=Y+\bar Y$ where $Y$ is a $(1,0)$ vector field, and let $u=0$ be a local equation for the divisor $C$. Then 
$$\omega+i\omega'=\frac{1}{u}\nu$$
where $\nu$ is a local non-vanishing holomorphic 2-form. Also, $\Vert \sigma \Vert^2= hu\bar u$ for a locally defined positive function $h$. 

If $i_X\omega=df$,
$$\frac{1}{u}i_Y\nu = \partial f= \partial \log hu\bar u =\frac{\partial h}{h}+\frac{du}{u}$$
and so
\begin{equation}
i_Y\nu= u\frac{\partial h}{h}+du
\label{Y}
\end{equation}
Since $\nu$ is nondegenerate this means that $Y$ (and hence $X$) is smooth in a neighbourhood of $C$. We can therefore integrate the globally defined vector field $X$ on the compact manifold $M$ to give a one-parameter group of diffeomorphisms $\varphi_t$. The form $\varphi_t^*(\omega+i\omega')=\omega+i\omega''$ is then meromorphic with respect to the complex structure $I^-=\varphi_t^*I$. 

We need to address next the behaviour of the form $(\varphi_t^*\omega'')^{1,1}$ -- does it extend to $M$, and is it positive?

Here we have on $U$
$$\frac{\partial }{\partial t}(\varphi_t^*\omega')^{1,1}=(\varphi_t^* {\mathcal L}_X\omega')^{1,1}=(\varphi_t^* dd^cf)^{1,1}.$$
But $dd^c\log \Vert \sigma \Vert^2=F$ where $F$ is the curvature of the connection on $K^*$ defined by the chosen hermitian metric. Thus $\omega^+(t)=(\varphi_t^*\omega')^{1,1}$ satisfies the differential equation
$$\frac{\partial \omega^+}{\partial t}=(\varphi_t^* F)^{1,1}$$
which is well-defined and smooth on the whole of $M$ since $\varphi_t$ and $F$ are globally defined. With the initial condition $\omega^+(0)=0$ this has the solution 
$$\omega^+(t)=\int_0^t(\varphi_t^* F)^{1,1}dt.$$
As before, if $F$ is a positive $(1,1)$ form then for small enough $t$, $\omega^+(t)$ will be positive and define a bihermitian metric. But there is a metric on $K^*$ with $F$ positive  if $K^*$ is ample, which is the case for a Del Pezzo surface.

\begin{rmk}
Surfaces which are not Del Pezzo can admit bihermitian structures -- indeed the quotient construction of \cite{T} yields Hirzebruch surfaces or $\CP^2$ blown up at an arbitrary number of points (though in special position). Our construction here has the property that the two complex structures $I^+,I^-$ are different but {\it equivalent} under a diffeomorphism.  This restricts the possibilities considerably. 

As an example consider $M$ to be the Hirzebruch surface $F_2$ -- this is the projective bundle $P({\mathcal O}\oplus {\mathcal O}(-2))\rightarrow \CP^1$. The canonical symplectic form on the cotangent bundle  ${\mathcal O}(-2)$ of $\CP^1$ extends as an anticanonical section with a double zero on the infinity section, which is a rational curve $C$ of self-intersection $+2$. Since $M\setminus C\cong T^*\CP^1$, the homology group $H_2(M\setminus C)$ is generated by the zero section $B$, a rational curve  of self-intersection $-2$. Now, for any bihermitian structure with this anticanonical divisor,  $\omega'$ and $\omega''$ are closed and smooth on $M\setminus C$. The integral of $\omega'$ on $B$ is zero because $\omega'$ is of type $(2,0)+(0,2)$ and $B$ is holomorphic, thus the cohomology class of $\omega'$ vanishes. On the other hand, from Proposition \ref{P2}, the integral of $\omega''$ is non-zero and so its cohomology class is non-zero. It follows that for the complex structure $I^-$, there can be no holomorphic curve in $M\setminus C$ and the complex structure must be different -- in fact it can only be $\CP^1\times \CP^1$. 
\end{rmk}

\section{Symplectic aspects}

The bihermitian 4-manifolds we have been considering are not obviously symplectic, but indeed they do have a natural symplectic structure. Note that $\omega'$ and $\omega''$ are closed  2-forms  which are singular along $C$. We shall see that the difference $\omega''-\omega'$ is a smooth symplectic form on $M$.

For the Del Pezzo surfaces, this is clear: the closed form $\rho(t)=\varphi_t^*\omega'-\omega'$ satisfies the equation
$$\frac{\partial \rho}{\partial t}=\varphi_t^* F$$
with initial condition $\rho(0)=0$ and so, by the same reasoning as above, is globally defined. Moreover $\rho^2=t^2 F^2+\dots$ and so for small $t$, since $F$ is positive, $\rho$ is a symplectic form, with cohomology class $2\pi tc_1(M)$. The result is true more generally:

\begin{prp} Let $(M,I)$ be a K\"ahler surface with a bihermitian structure such that $I^+=I$. Then
the 2-form $\omega''-\omega'$ extends to a smooth symplectic form $\rho$ on $M$. If $c_1(M)\ne 0$, the cohomology class of $\rho$ lies in the K\"ahler cone.
\end{prp}

  \begin{rmk} If $(M,I^+)$ is K\"ahler then $b_1(M)$ is even and from the classification in \cite{Ap} $(M,I^-)$ is K\"ahler, so the symplectic form lies in both K\"ahler cones.
\end{rmk}

\begin{prf}
From (\ref{fi1}) and (\ref{fi2}),
$$\omega''-\omega'=\frac{\phi'-\phi''}{8(1-p^2)}=\frac{\omega^++\omega^-}{4(1+p)}.$$
As shown in Proposition 2 of \cite{Ap}, when $b_1(M)$ is even (which holds for K\"ahler manifolds), $p\ne -1$, so that $\rho$ is well-defined on $M$ as a smooth closed 2-form. 

Now $(\omega^+)^2=(\omega^-)^2$ is the metric volume form and the $(1,1)$ component of $\omega^-$ is a multiple of $\omega^+$ because they are self-dual forms. From (\ref{fi2}) this multiple is $p$. Hence
$$\rho^2=\frac{(\omega^++\omega^-)^2}{16(1+p)^2}=\frac{(\omega^+)^2}{8(1+p)}$$
which is non-zero and hence $\rho$ is symplectic.

If $c_1(M)\ne 0$, then $\sigma$ vanishes somewhere, so there can be no holomorphic 2-forms and $h^{2,0}=0$, which means that the cohomology class $R=[\rho]$ is of type $(1,1)$. Moreover, from the classification of surfaces, $M$ is algebraic (see \cite{Ap}).

Restricting $\rho$ to a holomorphic curve with respect to $I^+$, $\omega'$ vanishes and $\omega''$ is positive from Proposition \ref{P2}. Hence $\rho$ is positive.  Since $R\cdot D>0$ for any effective divisor and $R^2>0$, by Nakai's criterion, the class $R$ is a K\"ahler class.
\end{prf}

The symplectic form $\rho$ defines a generalized complex structure on $M$, and (see \cite{H1}) it is, up to a B-field, one of the commuting pairs of generalized complex structures which define the bihermitian structure. The other one changes type from symplectic to complex on the elliptic curve $C$.

\section{Noncommutative geometry}

Suppose that the section $\sigma$ of $K^*$ vanishes with multiplicity one on a smooth elliptic curve $C$. Then  $C$ is given locally by a holomorphic function $u=0$ where $du\ne 0$ on $C$. From equation (\ref{Y}) we see that $X$ is the real part of a non-vanishing holomorphic vector field on $C$. Thus, the diffeomorphism $\varphi_t$ not only maps $C$ to itself, but does so by a translation. This means that, although $I^+$ and $I^-$ are equivalent complex structures on $M$, the restrictions of the  line bundles to $C$ differ by a translation. As an example, if we take $M$ to be $\CP^2$ then if a projective line with respect to  $I^+$ meets the cubic $C$ in the three points $x_1,x_2,x_3$, then using the addition law on the cubic curve, $x_1+ct,x_2+ct,x_3+ct$ are the points of intersection of a line in complex structure $I^-$.

This data -- an elliptic curve, a line bundle and a translation -- is the data for defining a noncommutative projective plane, as in \cite{Staff}. It would be interesting to know what role the non-holomorphic extension of this translation to $\CP^2$, which is the basis of our construction of bihermitian metrics, plays in noncommutative geometry from the differential geometric point of view.

If we think of $t\rightarrow 0$ as being the commutative limit, it is instructive to see what happens to the structures we have defined here. By definition, the complex structure $I^-(t)$ tends to $I^+$. As to the bihermitian metric $g(t)$, recall from Proposition \ref{P2} that the hermitian form is the $(1,1)$ part of $\varphi_t^*\omega'$. This tends to zero but
$$\lim_{t\rightarrow 0}\frac{1}{t} \varphi_t^*\omega'={\mathcal L}_X\omega'=F.$$
Thus if $g_0$ is the K\"ahler metric defined by $F$ and $I^+=I$, we have 
$$\lim_{t\rightarrow 0}\frac{1}{t} g(t)=g_0.$$

\vskip 1cm
 Mathematical Institute, 24-29 St Giles, Oxford OX1 3LB, UK
 
 hitchin@maths.ox.ac.uk
\end{document}